\documentclass[1p,authoryear]{elsarticle}

\usepackage{amsfonts}
\usepackage{amsmath}
\usepackage{amssymb}

\usepackage{abstract}

\usepackage{natbib}
\usepackage{geometry}

\usepackage{multicol}

\usepackage{graphicx}
\usepackage{subfigure,psfrag}

\newcommand{\vs}[1]{\vspace{#1}}

\newcommand{\begitem}{\begin{itemize}}    
\newcommand{\finit}{\end{itemize}}    
\newcommand{\begenum}{\begin{enumerate}}    
\newcommand{\finenum}{\end{enumerate}}    
\newcommand{\begar}{\begin{array}}    
\newcommand{\finar}{\end{array}}    
\newcommand{\begeq}[1]{\begin{equation} \label{#1}}    
\newcommand{\fineq}{\end{equation}}    
\newcommand{\begct}{\begin{center}}    
\newcommand{\finct}{\end{center}}    
\newcommand{\ra}{\rightarrow}    
\newcommand{\lra}{\longrightarrow}

\newcommand{\zun}{\vs{0.1cm}} 
     
\newcommand{\zdeux}{\vs{0.2cm}} 
\newcommand{\ztrois}{\vs{0.3cm}}    
\newcommand{\zcinq}{\vs{0.5cm}} 

\def\tinf{{\rightarrow\infty}}     

\def\bE{{\mathbb{E}}}

\def\bR{{\mathbb{R}}}

\def\bP{{\mathbb{P}}}

\def\tinf{{ \rightarrow \infty }}

\begin{document}
\setlength{\abovedisplayskip}{0.25cm}
\setlength{\belowdisplayskip}{0.25cm}

\newtheorem{theo}{Theorem}
\newtheorem{prop}{Proposition}
\newtheorem{defi}{Definition}
\newtheorem{lem}{Lemma}
\newtheorem{cor}{Corollary}
\newdefinition{rmk}{Remark}

\newcounter{Rem}   
\setcounter{Rem}{-1}    
\refstepcounter{Rem}    
\newcommand{\itemrem}{\addtocounter{Rem}{1} \underline{\rm Remarque \theRem} \ }   
\newcommand{\rem}{\addtocounter{Rem}{1}\hop\underline{\rm Remarque \theRem} \ }

\begin{center}
{\Large {\sc Empirical likelihood based confidence regions for}} \\
{\Large {\sc  first order parameters of heavy-tailed distributions } }
\bigskip

\large  Julien Worms(1) \& Rym Worms (2)

\medskip

 (1) Universit\'e de Versailles-Saint-Quantin-En-Yvelines\\
 Laboratoire de Math\'ematiques de Versailles (CNRS UMR 8100), \\
 UFR de Sciences, B\^at. Fermat, \\
 45 av. des Etats-Unis, 78035 Versailles Cedex, \\
 e-mail : {\tt worms@math.uvsq.fr}
\medskip

(2) Universit\'e Paris-Est-Cr\'eteil\\
Laboratoire d'Analyse et de Math\'ematiques Appliqu\'ees (CNRS UMR 8050), \\
 61 av. du G\'en\'eral de Gaulle, 94010 Cr\'eteil cedex, \\
 e-mail : {\tt rym.worms@u-pec.fr}
\end{center}

  \vspace{0.1cm} 
\setlength{\belowdisplayskip}{0.08cm}
\setlength{\abovedisplayskip}{0.cm}

\begin{abstract} 
Let $X_1, \ldots, X_n$ be some i.i.d. observations from a heavy tailed distribution $F$, i.e. such that the common distribution of the excesses over a high threshold $u_n$ can be approximated by a Generalized Pareto  Distribution $G_{\gamma,\sigma_n}$ with $\gamma >0$. This paper deals with the problem of finding confidence regions for the couple $(\gamma,\sigma_n)$ : combining the empirical likelihood methodology with estimation equations (close but not identical to the likelihood equations) introduced by \cite{Zhang07}, asymptotically valid confidence regions for $(\gamma,\sigma_n)$ are obtained and proved to perform better than Wald-type confidence regions (especially those derived from the asymptotic normality of the maximum likelihood estimators). By profiling out the scale parameter, confidence intervals for the tail index are also derived. 
\end{abstract}

\setlength{\belowdisplayskip}{0.19cm}
\setlength{\abovedisplayskip}{0.19cm}

\section{Introduction}
 \label{intro}
 
 \noindent In statistical extreme value theory, one is  often interested by the estimation of the so-called tail index $\gamma=\gamma(F)$ of the 
underlying model $F$ of some i.i.d. sample $(X_1,\ldots, X_n)$, which is the shape parameter of the Generalized Pareto Distribution (GPD) with distribution function (d.f.)
\[ 
 G_{\gamma, \sigma} (x) = \left\{ \begar[c]{ll}
 1-\Bigr(1 \, + \, {\gamma x \over \sigma}\Bigr)^{-{1 \over \gamma}}, & \mbox{for } \gamma \neq 0 \\
 1-\exp\Bigr( - {x  \over \sigma} \Bigr), & \mbox{for } \gamma=0. \finar \right.
\]
 
\noindent The GPD appears as the limiting d.f. of excesses over a high threshold $u$ defined for $x\geq 0$ by  
 \[
 F_u(x):=\mathbb{P}(X-u\leq x \,\vert\, X>u), \mbox{ \ where $X$ has d.f. $F$. } 
 \]

\noindent
It was established in \cite{Pickands75}  and  \cite{BalkemaHaan74}
that $F$ is in the domain of attraction of an extreme value distribution with shape parameter $\gamma$  if and only if
\begeq{theoPickands}
 \lim_{u \to s_+(F)} \sup_{0<x<s_+(F)-u} \Bigr\vert F_u(x) - G_{\gamma, \sigma(u)} (x)\Bigr\vert = 0 
\fineq

\noindent
for some positive scaling function $\sigma(\cdot)$, where $s_+(F)=\sup \{ x: F(x)<1 \}$. 
This suggests to model the  d.f. of excesses over a high threshold by a GPD. This is the P.O.T method. 
\zdeux

\noindent
Some estimation methods for the couple $(\gamma, \sigma)$ in the GPD parametrization have been proposed. We can cite the 
maximum likelihood (ML) estimators of  \cite{Smith87} or the probability weighted moments (PWM) estimators of 
\cite{HoskWall87}. In \cite{Zhang07}, the author proposed new estimators based on estimating equations close to the likelihood equations. 
Using the reparametrization $b=-\gamma/\sigma$ and considering $X_1,\ldots,X_n$ i.i.d. variables with distribution $G_{\gamma,\sigma}$ with $\sigma$ a fixed value (which is an important restriction if the aim is to prove asymptotic results), 
he based his method on one of  the likelihood equations 
\begeq{equMLE}
\gamma= -\frac{1}{n} \sum_{i=1}^n  \log (1-bX_i)
\fineq 
and on the empirical version of the moment equation $\bE((1-bX_1)^r)=\frac{1}{1-r \gamma}$, i.e. 
\[ 
\frac{1}{n} \sum_{i=1}^n (1-bX_i)^r -\frac{1}{1-r \gamma} =0
\]
or 
\begeq{equMom}
\frac{1}{n} \sum_{i=1}^n (1-bX_i)^{r/\gamma} -\frac{1}{1-r } =0,  
\fineq 
for some parameter $r < 1$. $(\ref{equMLE})$ and $(\ref{equMom})$ yield the estimation equation for $b$ 
\begeq{equb}
\frac{1}{n} \sum_{i=1}^n (1-bX_i)^{nr (\sum_{i=1}^n  \log (1-bX_i))^{-1}} -\frac{1}{1-r } =0, \ \ \mbox{ provided } b < X_{(n)}^{-1} \mbox{ and } \ r < 1. 
\fineq
An estimation of $\gamma$ is then deduced from $(\ref{equMLE})$ and $\sigma$ is estimated using $b=-\gamma/\sigma$ .
\zun\\
Zhang proved in \cite{Zhang07}  (Theorem 2.2) that for a GPD$(\gamma,\sigma)$ sample with $\gamma > -1/2$, the estimators he proposed for $\gamma$ and $\sigma$ are  jointly asymptotically normally distributed and that they share none of the following drawbacks of the ML and PWM methods : theoretical invalidity of large sample results for the PWM estimators with large positive $\gamma$, 
and computational problems for the ML estimator.  
\ztrois\\
In this paper, we consider the classical P.O.T. framework, where an i.i.d. sample $X_1, \ldots, X_n$ with distribution $F$ is observed and, according to ($\ref{theoPickands}$), a GPD $G_{\gamma,\sigma(u_n)}$ is fitted to the sample of the excesses over a large threshold $u_n$. Noting $\sigma_n= \sigma(u_n)$, our goal is to build confidence regions for the couple $(\gamma,\sigma_n)$ (as well as  confidence intervals for the tail index $\gamma$)
for heavy-tailed distributions ($\gamma >0$ case), starting from Zhang's estimating equations; therefore the excesses will be approximately GPD distributed and the parameter $\sigma=\sigma_n$ will be varying with $n$. To the best of our knowledge, little attention has been paid to the subject of joint confidence regions (and their coverage probabilities) for the couple $(\gamma,\sigma_n)$, especially outside the exact GPD
framework. 
\zun\\
An obvious approach to obtain confidence regions is to use the gaussian  approximation. In this work, we consider an alternative method, namely the empirical likelihood method. This method was developped by Owen in \cite{Owen88} and \cite{Owen90}, for the mean vector of i.i.d. observations and has been extended to a wide range of applications, particularly for generalized estimating equations (\cite{QinLawless94}) . 
\zun\\
In \cite{LuPeng02}, this method was applied to construct  confidence intervals for the tail index of a heavy-tailed distribution (empirical likelihood estimator of the tail index being equal to the Hill estimator). It turned out that the empirical likelihood method performs better than the normal approximation method in terms of coverage probabilities especially if the calibration method proposed in \cite{PengQi06} is adopted.  We will see that it is even more the case for confidence regions. 
\zun\\
In Section 2, we explain the empirical likelihood methodology based on Zhang's equations $(\ref{equMLE})$ and $(\ref{equb})$ , and present some asymptotic results.  A simulation study is conducted in Section 3, which compares different methods for constructing confidence regions for the couple $(\gamma, \sigma_n)$, as well as  confidence intervals for $\gamma$ alone, in terms of coverage probabilities. Proofs are given in Section 4 and some details left in the Appendix. Technical difficulties are mainly due to the fact that one of the parameters and the distribution of the excesses are depending  on $n$. 

\section{Methodology and statement of the results}

\subsection{Notations and Assumptions}
\label{hypotheses}

In this work, the tail index $\gamma$ is supposed to be positive and $F$  twice differentiable with well defined inverse $F^{-1}$. 
Let $V$ and $A$ be the two functions defined by
\[
 U(t)=\bar{F}^{-1}(1/t) \quad \hbox{and} \quad A(t) \,  = t \frac{U^{\prime \prime}(t)}{U^{\prime}(t)} +1 - \gamma,
\]
where $\bar{F}=1-F$.
 
\vskip0ex
\noindent
We suppose the following first and second order conditions hold ($RV_{\rho}$ below stands for the set of  regularly varying functions with coefficient  of variation $\rho$) :
\begin{eqnarray} 
& \lim_{t \rightarrow +\infty} A(t) = 0 & \label{cond1erAf} \\ 
& A \; \hbox{ is of constant sign at $\infty$ and there exists } \;
\rho \leq 0 \; \hbox{ such that } \;  \vert A \vert \in RV_{\rho}. & \label{cond2deARho} 
\end{eqnarray}

\noindent
A proof of the  following lemma can be found in \cite{Haan84}. 
\begin{lem}
\label{resultatsfoncU}
Under $(\ref{cond1erAf})$ and $(\ref{cond2deARho})$ we have, for all $x>0$, 
\begeq{cond2U}
 \left( \frac{U(tx)-U(t)}{tU^{\prime} (t)}-\frac{x^{\gamma}-1}{\gamma} \right) \big/ A(t)  \longrightarrow   K_{\gamma, \rho}(x), 
 \mbox{ as } t \rightarrow +\infty, 
\fineq 
where $K_{\gamma, \rho}(x):= \int_1^x  u^{\gamma-1} \int_1^u s^{\rho-1} ds du$,
 and the following well-known Potter-type bounds hold: \zun 
 
 \noindent 
 $\forall \epsilon >0,  \  \exists t_0, \ \forall t \geq t_0, \  \forall x  \geq 1$,
\begeq{bornesPotter}
  (1-\epsilon) \exp^{-\epsilon \log(x)} K_{\gamma, \rho}(x)  
   \; \leq \;  \left( \frac{U(tx)-U(t)}{tU^{\prime} (t)}-\frac{x^{\gamma}-1}{\gamma} \right) \big/ A(t) 
   \; \leq \; (1+\epsilon) \exp^{\epsilon \log(x)} K_{\gamma, \rho}(x).
\fineq 
\end{lem}

\subsection{Confidence regions for the couple $(\gamma,\sigma_n)$}
\label{couple}
\zun 
For some $r<1$ and positive $y,\gamma,\sigma$, let 
\[
g(y,\gamma,\sigma) :=  \left( \begar{cc}    \log(1+\gamma  y / \sigma)-\gamma \zun \\ 
                                                        (1+\gamma y / \sigma)^{r/\gamma}-\frac{1}{1-r} \finar \right).
\]
Note that, if $Z_1, \ldots, Z_n$ are i.i.d. GPD$(\gamma,\sigma)$, then $\frac{1}{n} \sum_{i=1}^n g(Z_i,\gamma,\sigma)=0$ summarizes  equations $(\ref{equMLE})$ and $(\ref{equMom})$ of \cite{Zhang07}.  
\zun\\ 
Let $X_1, \ldots, X_n$ be i.i.d. random variables with common  d.f. $F$ (satisfying the assumptions stated in the previous paragraph),
and $\gamma_0$ and $\sigma_0(\cdot)$ be the true parameters such that relation (\ref{theoPickands}) is satisfied.  
For a fixed high threshold $u_n$, consider the $N_n$ excesses $Y_1\ldots, Y_{N_n}$ over $u_n$. Conditionally  on $N_n=k_n$,  
$Y_1\ldots, Y_{k_n}$ are i.i.d. with common distribution function $F_{u_n}$ which , according to  ($\ref{theoPickands}$), is  approximately 
$G_{\gamma_0,\sigma_{0n}}$, where $\sigma_{0n}:=\sigma_0(u_n)$. The objective is to estimate $\gamma_0$ and $\sigma_{0n}$.
\zun\\
Let ${\cal S}_n$ denote the set of all probability vectors $p=(p_1, \ldots, p_{k_n})$ such that $\sum_{i=1}^{k_n} p_i=1$ and $p_i \geq 0$.
The empirical likelihood for $(\gamma,\sigma)$ is defined by 
 \[
 L(\gamma,\sigma) := \sup \left\{ \, \prod_{i=1}^{k_n} p_i \, \left/ \; p \in {\cal S}_n \, \mbox{ and }  \, \sum_{i=1}^{k_n} p_i g(Y_i,\gamma,\sigma) =0 
 \, \right. \right\}
 \] 
and the empirical log likelihood ratio is then defined as 
\[ 
l(\gamma,\sigma) 
:= -2 (\log L(\gamma,\sigma) - \log L(\hat{\gamma}_n,\hat{\sigma}_n))
\] 
where $(\hat{\gamma}_n,\hat{\sigma}_n)$ are maximising $L(\gamma,\sigma)$, and are called the maximum empirical likelihood estimates (MELE)
of the true  parameters $\gamma_0$ and $\sigma_{0n}$. 
\zun\\
Since Theorem 2.1 of \cite{Zhang07} implies that, for $r<1/2$, there exists a unique and easily computable solution $(\tilde\gamma_n,\tilde\sigma_n)$ to the equations
\[
 \frac{1}{k_n} \sum_{i=1}^{k_n} \log \left( 1+\gamma Y_i/\sigma \right) - \gamma =0  \  \  
 \mbox{ and }Ê
 \frac{1}{k_n} \sum_{i=1}^{k_n}  \left( 1+\gamma Y_i/\sigma \right)^{r/\gamma} - \frac{1}{1-r}=0
\]
{\it i.e.} such that $k_n^{-1} \sum_{i=1}^{k_n} g(Y_i,\tilde\gamma_n,\tilde\sigma_n) = 0$, it thus comes that $L(\tilde\gamma_n,\tilde\sigma_n)
=k_n^{-k_n}$ which is equal to $\max_{\gamma,\sigma} L(\gamma,\sigma)$ : the MELE estimators $(\hat{\gamma}_n,\hat{\sigma}_n)$ therefore 
coincides with Zhang's estimators. 
\zun

\noindent   Note however that Zhang worked in the purely GPD framework and that the application of his results for 
constructing confidence regions, based on the asymptotic normality of $(\hat{\gamma}_n,\hat{\sigma}_n)$, necessarily involves 
some additional covariance estimation. 
Our aim is to construct confidence regions for $(\gamma,\sigma_n)$ directly, relying on the asymptotic distribution
of the empirical likelihood ratio $l(\gamma_0,\sigma_{0n})$ stated in the following theorem. Classical advantages of proceeding so are well known : a first one is the avoidance of information matrix estimation, a second one is the guarantee of having the confidence region included in the parameter space (as a matter of fact, in our case, the parameter $\sigma$ is positive but nothing guarantees that the confidence region for $(\gamma,\sigma)$, based on the CLT for $(\hat{\gamma}_n,\hat{\sigma}_n)$ will not contain negative values of $\sigma$) . Note in addition that our result is proved in 
the general framework ({\it i.e.} when the excess distribution function is supposed to be only approximately GPD) and that simulation results
show some improvements in terms of coverage probability (see next Section). 
\zdeux\\
Note that the empirical log-likelihood ratio $l(\gamma,\sigma) = -2 \log (k_n^{k_n} L(\gamma,\sigma))$ has a more explicit expression : following \cite{Owen90},  the Lagrange multipliers method yields
 \[
p_i=\frac{1}{k_n (1+<\lambda(\gamma,\sigma),g(Y_i,\gamma,\sigma)>)}  \ \  \mbox{ and }  \  \   l(\gamma,\sigma) = 2 \sum_{i=1}^{k_n} \log \left( 1+ <\lambda(\gamma,\sigma),g(Y_i,\gamma,\sigma)> \right),
 \]
where $\lambda(\gamma,\sigma)$ is determined as the solution of the system
 \begeq{equlambda}
 \frac{1}{k_n} \sum_{i=1}^{k_n} \left(  1+<\lambda(\gamma,\sigma),g(Y_i,\gamma,\sigma)>\right)^{-1} g(Y_i,\gamma,\sigma) =0.
 \fineq
Let  $a_n  := A \left( 1/\bar{F}(u_n) \right)$.
 
\begin{theo}
\label{asymptCouple}
Under conditions (\ref{cond1erAf}) and (\ref{cond2deARho}), with $\gamma >0$, conditionally on $N_n=k_n$, if we suppose that  $k_n$ tends to $+\infty$ such that $\sqrt{k_n} a_n$   goes to $0$  as $n \rightarrow +\infty$, then for $r<1/2$ 
\[
l(\gamma_0,\sigma_{0n}) \stackrel{\cal{L}}{\rightarrow} \chi ^2(2), \mbox{ as }  n \rightarrow +  \infty. 
\]
\end{theo}

\noindent This result is the basis for the construction of a confidence region, of asymptotic level $1-\alpha$, for the couple $(\gamma_0,\sigma_{0n})$ which consists in  all $(\gamma,\sigma)$ values such that $l(\gamma,\sigma) \leq c_\alpha$, where $c_\alpha$ is the $1-\alpha$ quantile  of the $\chi ^2(2)$ distribution. \zun \\ 
Note that $\sqrt{k_n} a_n \rightarrow 0$ was also assumed in \cite{LuPeng02}.

\subsection{Confidence interval for  $\gamma$}
\label{Idcgamma}

For a fixed parameter $\gamma$, we note $\hat{\sigma}_{\gamma}$ the value of $\sigma$ that minimizes $l(\gamma, \sigma)$. 
Then, $l(\gamma,\hat{\sigma}_{\gamma})$ is  called the profile empirical log likelihood ratio. The following asymptotical result is the basis for constructing the confidence intervals  for the true parameter $\gamma_0$ of the model. 
\begin{theo}
\label{asymptprofile}
Under the same conditions as Theorem $\ref{asymptCouple}$, if $r<1/3$ then, conditionnally on $N_n=k_n$, 
\[
 l(\gamma_0,\hat{\sigma}_{\gamma_0}) \stackrel{\cal{L}}{\rightarrow} \chi ^2(1), \mbox{ as }  n \rightarrow +  \infty. 
\]
\end{theo}

\noindent This result yields as a confidence interval with asymptotic level  $1-\alpha$ for the tail index $\gamma_0$, the set of  all $\gamma$ values such that $l(\gamma,\hat{\sigma}_{\gamma}) \leq c_\alpha$, where $c_\alpha$ is the $1-\alpha$ quantile  of the $\chi ^2(1)$ distribution.

\begin{rmk}
Note that the restriction $r< 1/3$ could be reduced to $r< 1/2$, but this would unnecessarily complicate the proof since most of the time $r$ should be chosen negative (see \cite{Zhang07} for a discussion). 
\end{rmk}

\section{Simulations}

\subsection{Simulations for the couple $(\gamma, \sigma)$}
\label{Simulcouple}

In this subsection, we present a small simulation study in order to investigate the performance of our proposed method  for constructing confidence regions for the couple $(\gamma_0, \sigma_{0n})$ based on empirical likelihood techniques (Theorem \ref{asymptCouple}).  We compare empirical coverage probabilities of the confidence regions (with nominal level $0.95$) produced by our empirical likelihood method (EL$(0.95)$), the normal approximation for the Maximum Likelihood  estimators (ML$(0.95)$) and the normal approximation for the estimators proposed in \cite{Zhang07} (Zhang$(0.95)$). 
\zdeux

\noindent
Before giving more details on the simulation study, let us make the following remarks.
\begin{rmk} The CLT for the GPD parameters stated in \citet{Zhang07}  has been proved  when the underlying distribution is a pure GPD : using Theorem 2.1 in \cite{Zhang07} (which asserts the existence and unicity of the estimator), and following the same methodology that led to 
Proposition \ref{proplocalminimizer} in Section \ref{preuve2}, the consistency of the sequence $(\hat\gamma,\hat\sigma_n)$ can be obtained 
(details are omitted here) and accordingly, by classical methods, the asymptotic normality follows in the general case of the underlying distribution 
belonging to the Fr\'echet maximum domain of attraction (assuming that $\sqrt{k_n}a_n\rightarrow 0$). We will therefore use the convergence in distribution of $\sqrt{k_n} ( \, \hat\gamma - \gamma_0 \, , \, \hat\sigma_n/\sigma_{0n} -1 \, )$ to ${\cal N}_2(0,\Sigma(\gamma_0,r))$, where  
\[
 \Sigma(\gamma,r) := \left( \begar[c]{cc} (1-r)(1+(2\gamma^2+2\gamma+r)/(1-2r)  &  -1 - (r^2+\gamma^2+\gamma)/(1-2r) \\
 -1 - (r^2+\gamma^2+\gamma)/(1-2r)  &  2 + ((r-\gamma)^2+2\gamma)/(1-2r)  \finar \right)  
\]
and consider the corresponding  Wald-type confidence region,  based on approximating the  distribution of  the statistic
$\|\sqrt{k_n} (\Sigma(\hat\gamma,r))^{-1/2}( \, \hat\gamma - \gamma_0 \, , \, \hat\sigma_n/\sigma_{0n} -1 \, )\|^2$ by $\chi^2(2)$. The same methodology is
applied for constructing confidence regions based on the Maximum Likelihood estimators (see \cite{Smith87}).
\end{rmk}
\begin{rmk} The tuning parameter $r$ in  $(\ref{equb})$ is chosen equal to $-1/2$ (as suggested in \cite{Zhang07}, without any prior information).
The empirical likelihood confidence region is based on a Fisher calibration rather than a $\chi^2$ calibration, as suggested in \cite{Owen90} :
concretely, this means that our confidence region consists in  the set of all $(\gamma,\sigma)$ such that $l(\gamma,\sigma)\leq f$ where
$\bP(\, (2(k_n-1)/(k_n-2)) F(2,k_n-2) \leq f \,) = 0.95$. Indeed, it has been empirically observed that it (uniformly in $k_n$) produces slightly better 
results in terms of coverage probability. 
\end{rmk}

\noindent The simulations below are based on 2000 random samples of size $n=1000$, generated from the following two distributions: the Fr\'echet  
distribution with parameter $\gamma >0$ given by $F(x)=\exp(-x^{-1/\gamma})$ (the results for $\gamma=1$ and $1/4$ are presented) 
and the Burr distribution with positive  parameters $\lambda$ and $\tau$ given by $F(x)=(1+x^{\tau})^{-1/\lambda}$ 
(the results for $(\lambda,\tau)=(1,1)$ and  $(2,2)$ are presented). Note that the tail index for the Burr distribution 
is $\gamma= (\lambda  \tau)^{-1}$.  
\zdeux

\begin{figure}[htp]
    \centering
     \subfigure[Coverage Probability for Burr$(1,1)$ model, $n=1000$]{
          \label{Burr11RdC}
          \includegraphics[width=.45\textwidth]{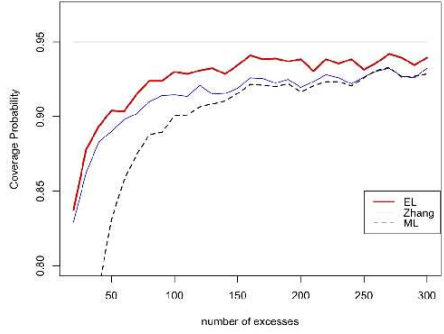}}
     \hspace{.1in}
     \subfigure[Coverage Probability for Burr$(2,2)$ model, $n=1000$]{
          \label{Burr22RdC}
          \includegraphics[width=.45\textwidth]{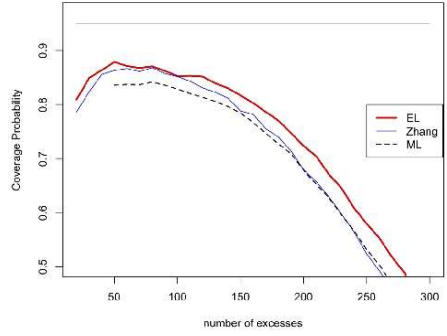}}
     \subfigure[Coverage Probability for Frechet$(1)$ model, $n=1000$]{
          \label{Frechet1RdC}
          \includegraphics[width=.45\textwidth]{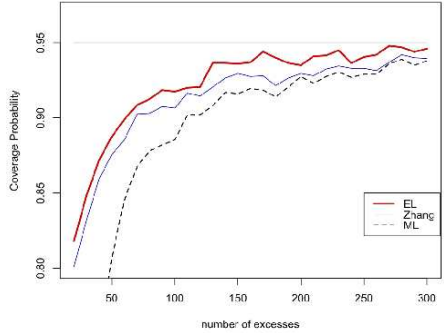}}
     \hspace{.1in}
     \subfigure[Coverage Probability for Frechet$(1/4)$ model, $n=1000$]{
          \label{Frechet1quartRdC}
          \includegraphics[height=6.8cm,width=.45\textwidth]{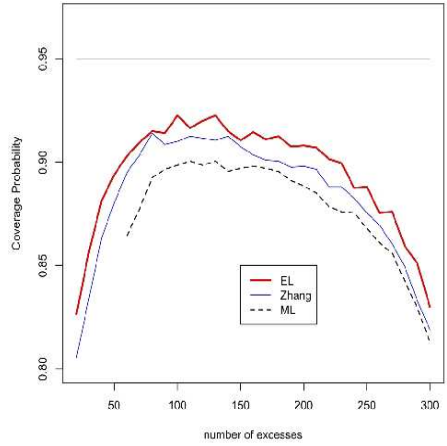}}
     \caption{Coverage Probability for Burr$(1,1)$, Burr$(2,2)$, Frechet$(1)$, Frechet$(1/4)$ as a function of the number of excesses $k_n$.
     The dashed line is for ML, the thin solid line for Zhang, and the thick solid line for EL.}
     \label{figure1}
\end{figure}

\noindent Coverage probabilities EL$(0.95)$, ML$(0.95)$ and Zhang$(0.95)$ are plotted against different values of $k_n$, the number 
of excesses used.  Figure \ref{figure1} seems to indicate that our method performs  better in terms of coverage probability, and additional
models not presented here almost always led to the same ordering : EL better than the bivariate CLT of 
Zhang's estimator, itself better than the CLT for the MLE. However, we have observed that the overall results are not very satisfactory 
when the tail index $\gamma$ is small.
\zdeux

\begin{figure}
 \centering
\includegraphics[width=.45\textwidth]{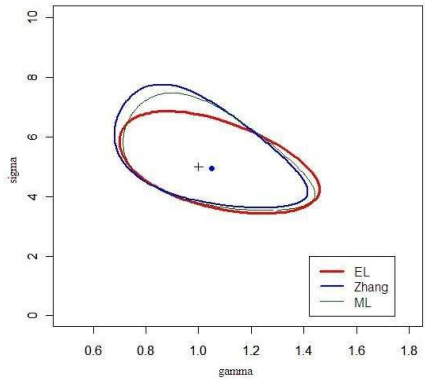}
 \caption{Confidence regions for  a sample of Burr$(1,1)$. The thinner line  is for ML, the thicker one is  for EL and the other  for Zhang. }

   \label{figure2}
 \end{figure}

\noindent One can wonder if the improvement  in coverage probabilities is due to the fact that our confidence regions are wider than in the ML and Zhang cases. In practice, it seems in fact that the three confidence regions have comparable sizes (our confidence region being even a bit smaller). Figure \ref{figure2} shows these three regions for a simulated Burr$(1,1)$ random sample with $n=1000$ et $k_n=200$.

\begin{rmk} it should be noted that some computational problems occurred when trying to calculate the maximum likelihood estimators. This explains why in some of the figures above, a piece of the curve is lacking for the smaller values of $k_n$ (the computation was performed
by using the function {\tt fpot} of the {\tt evd} package of R). 
\end{rmk}

\subsection{Simulations for $\gamma$}
\label{Simulgamma}

In this subsection, we present another small simulation study which is now concerned by the performance of our method for constructing confidence intervals  for the tail index  $\gamma$ based on profile empirical likelihood techniques (Theorem \ref{asymptprofile}).  We compare empirical coverage probabilities of the confidence intervals  for our empirical  profile likelihood method (ELW$(0.95)$), the empirical likelihood method proposed in \cite{LuPeng02} and based on the Hill estimator (ELP$(0.95)$) and finally the normal approximation for the estimator of $\gamma$ proposed in \cite{Zhang07} (Zhang$(0.95)$). Note 
that for ELP$(0.95)$, we used the exponential calibration in order to calculate the critical value, as prescribed in \cite{PengQi06}.
As before, the Fisher calibration  was preferred to the $\chi^2$ one in order to compute the critical value in our case. 
We worked with the same two distributions as in the couple case : the results are presented for $\gamma=1$ and $1/4$ in the Fr\'echet case and $(\lambda,\tau)= (1,1)$ and $(2,2)$  in the Burr one.
\zdeux

\begin{figure}
    \centering
    \subfigure[Coverage Probability for Burr$(1,1)$ model, $n=1000$]{
          \label{Burr11IdC}
          \includegraphics[height=7.0cm,width=.45\textwidth]{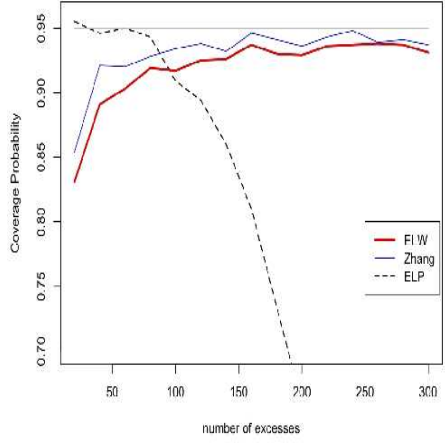}}
     \hspace{.1in}
     \subfigure[Coverage Probability for Burr$(2,2)$ model, $n=1000$]{
          \label{Burr21IdC}
          \includegraphics[width=.45\textwidth]{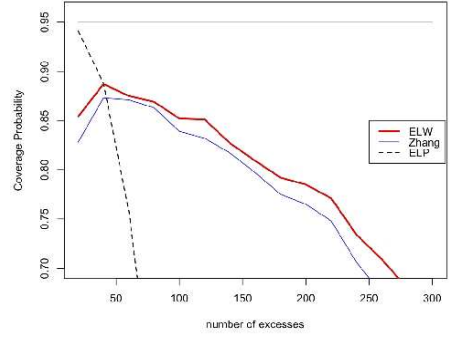}}
     \subfigure[Coverage Probability for Frechet$(1)$ model, $n=1000$]{
         \label{Frechet1IdC}
        \includegraphics[height=6.8cm,width=.45\textwidth]{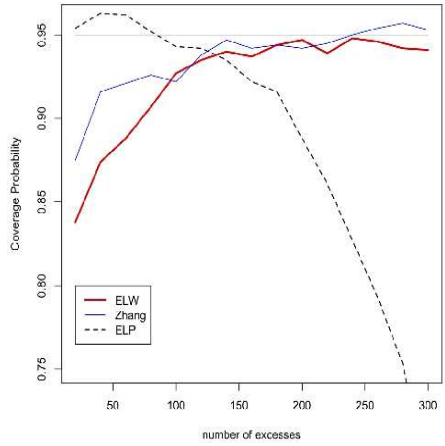}}
     \hspace{.1in}
    \subfigure[Coverage Probability for Frechet$(1/4)$ model, $n=1000$]{
         \label{Frechet1quartIdC}
          \includegraphics[height=7.0cm,width=.45\textwidth]{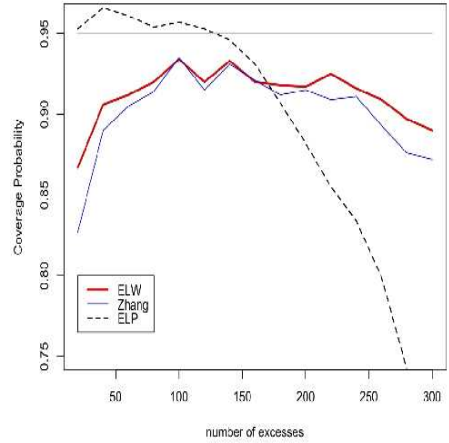}}
     \caption{Coverage Probability for Burr$(1,1)$, Burr$(2,2)$, Frechet$(1)$, Frechet$(1/4)$ as a function of the number of excesses $k_n$.
     The dashed line is for $ELP$, the thin solid line for $Zhang$, and the thick solid line for $ELW$.}
     \label{figure3}
\end{figure}

\noindent  Empirical coverage probabilities are plotted against different values of $k_n$. Contrary to the couple framework, Figure \ref{figure3}  shows no  significant improvement with respect to Zhang's CLT based confidence intervals, which itself 
shows some problems of undercoverage. The EL-based confidence intervals of Lu, Peng and Qi, show quite satisfactory coverage probabilities, but on a range of values of $k_n$ which is sometimes very narrow (which is a common phenomenon in POT methodology) : this drawback is much
less present for the confidence intervals based on Zhang's equations, which show a better stability against the choice of $k_n$. Moreover, simulations
(not presented here) showed that the widths of the ELP interval and the ELW interval are comparable, and smaller than the width of the interval
based on Zhang's CLT.  
\zdeux

\begin{rmk} the computation of the profile empirical likelihood $l(\gamma_0,\hat\sigma_{\gamma_0})$ was performed using a classical 
descent algorithm, taking profit of some convexity properties of the profile empirical likelihood function. Computational details and files can be 
obtained from the authors (some of them are downloadable on the first author's webpage).
\end{rmk}

\section{Proofs}

Note that we will prove Theorem \ref{asymptprofile} before Theorem \ref{asymptCouple} because its proof is more involved 
and largely includes what is needed to prove Theorem \ref{asymptCouple}.

\subsection{Proof of Theorem \ref{asymptprofile}}
\label{preuve2}
From now on we work conditionally on $\{N_n=k_n\}$ for some given sequence $(k_n)$ satisfying $\sqrt{k_n}a_n\ra 0$ as $n\tinf$.  
\zun\\
Let $U_1,\ldots,U_{k_n}$ denote independent uniform random variables on $[0,1]$. Noticing that $(Y_1,\ldots,Y_{k_n})$ has the
same joint distribution as $(\tilde Y_1,\ldots,\tilde Y_{k_n})$ defined by 
\[
 \tilde Y_i := \bar{F}^{-1}_{u_n}(U_i) = U(1/(U_i \bar{F}(u_n)))-U(1/\bar{F}(u_n))
\]
we see that it suffices to prove Theorem \ref{asymptprofile} with $\tilde Y_i$ replacing $Y_i$ in the definition of the empirical
likelihood. {\it For simplicity we will write $Y_i$ instead of $\tilde Y_i$ in the sequel} (note that these random variables are
i.i.d. but their common distribution varies with $n$).
\zun\\
Now, defining $Z_{i,n}:=Y_i/\sigma_{0n}$, $\lambda_0(\theta):=\lambda(\gamma_0,\sigma_{0n}\theta)$ 
and $g_0(z,t):=g(z,\gamma_0,t)$, we see that $g_0(Z_{i,n},\theta) = g(Y_i,\gamma_0,\sigma_{0n}\theta)$, hence
\[ \begar{ll}
 l_0(\theta) \, := \, l(\gamma_0,\sigma_{0n}\theta)
 &  \, = \, 2\sum_{i=1}^{k_n} \log \left( 1+<\lambda(\gamma_0,\sigma_{0n}\theta),g(Y_i,\gamma_0,\sigma_{0n}\theta)> \right) \zun \\
&  \, = \, 2\sum_{i=1}^{k_n} \log \left( 1+<\lambda_0(\theta),g_0(Z_{i,n},\theta)> \right). 
\finar \]
{\it
With these preliminaries in mind, we thus need to prove that there exists some local minimizer $\hat\theta$ of $l_0(\cdot)$ in a neighborhood
of $\theta_0=1$ such that $l_0(\hat\theta)\ra \chi^2(1)$ in distribution, because $l_0(\hat\theta)=l(\gamma_0,\hat\sigma_{\gamma_0})$ with 
$\hat\theta=\hat\sigma_{\gamma_0}/\sigma_{0n}$.
} 
\zun\\
We now state in the following proposition some important results, which will be proved in Section \ref{preuveprop} and will enable us to proceed with 
the proof of Theorem \ref{asymptprofile} following a plan very similar to that found in \cite{QinLawless94} 
(note that here the parameter is one-dimensional whereas the estimating function $g_0$ is $\bR^2$-valued). We first introduce some
important notations :
\[
 \begar{ccl}
  G_n(\theta)& := & \frac{1}{k_n} \sum_{i=1}^{k_n} g_0(Z_{i,n},\theta) \zun \\
  B_n(\theta) & := & \frac{1}{k_n} \sum_{i=1}^{k_n} g_0(Z_{i,n},\theta)(g_0(Z_{i,n},\theta))^t \zun\\
  A_n(\theta) & := & \frac{1}{k_n} \sum_{i=1}^{k_n} \frac{\partial g_0}{\partial\theta} (Z_{i,n},\theta) \zun\\
  M_n(\theta) & := & \max_{i\leq k_n} \|g_0(Z_{i,n},\theta)\|
 \finar  
\]  
Note that, although the new parameter $\theta$ is scalar, we will write below $\|\theta\|$ instead of $|\theta|$
in order to emphasize the fact that the arguments described below can be applied to more general frameworks. The same is true about
the fact that we use below the notation  $\theta_0$ instead of simply the number $1$.

\begin{prop}
\label{propprelim}
Suppose the assumptions of Theorem \ref{asymptprofile} are valid, and $Z$ is some random variable
distributed as $G_{\gamma_0,1}$. If ${\cal B}_n:=\{ \, \theta\in\bR \, ; \, \|\theta-\theta_0\| \leq k_n^{-1/3} \, \}$, 
then we have, conditionally to $\{N_n=k_n\}$ and as $n\tinf$, 
\zun
\begin{eqnarray}
 G_n(\theta_0) & = & O\left( k_n^{-1/2}(\log\log k_n)^{1/2} \right) \mbox{ a.s.} 
 \label{Gntheta0vitesseps} 
 \zun\\
 \sqrt{k_n} G_n(\theta_0) & \stackrel{d}{\longrightarrow} & {\cal N}(0,B) 
 \label{Gntheta0TLC}
 \\
 M_n(\theta_0) & = & o(\sqrt{k_n}) \mbox{\ a.s.  }
 \label{controleMntheta0}
 \ztrois\\
 M_n(\theta) & = & o(k_n^{1/3}) \mbox{\ a.s. uniformly on ${\cal B}_n$}
 \label{controleMntheta}
 \\
 A_n(\theta) & = & A + o(1) \mbox{\ a.s. uniformly on ${\cal B}_n$}, 
 \label{cvgceunifAntheta}
 \zun\\
 B_n(\theta) & = & B + o(1) \mbox{\ a.s. uniformly on ${\cal B}_n$}
 \label{cvgceunifBntheta}
 \ztrois\\
 G_n(\theta) & = & G_n(\theta_0)+(A+o(1))(\theta-\theta_0) = O(k_n^{-1/3}) \mbox{\ a.s. unif. on ${\cal B}_n$}
 \label{controleGntheta}
 \end{eqnarray} \vspace{-0.1cm}
 where 
\[ B \; := \; \bE\left[ \,  g_0(Z,\theta_0)(g_0(Z,\theta_0))^t \, \right] \; = \; 
 \left( \begar{cc}  \gamma^2 & \frac{\gamma r}{(1-r)^2} \\ \frac{\gamma r}{(1-r)^2} & \frac{r^2}{(1-2r)(1-r)^2} \finar \right)
 \]
 and 
 \[
 A \; := \; \bE\left[ \,  \frac{\partial g_0}{\partial\theta} (Z,\theta_0) \, \right] \; = \; 
 \left( \begar{c}  -\frac{\gamma}{\gamma+1}  \\  -\frac{r}{(1-r+\gamma)(1-r)}  \finar \right)
 \]
\end{prop}

Recall that $\lambda_0(\theta)$ is defined as the solution of the equation 
\[
  \frac{1}{k_n} \sum_{i=1}^{k_n} \left(  1+<\lambda,g_0(Z_{i,n},\theta)>\right)^{-1} g_0(Z_{i,n},\theta) = 0.
\] 
Therefore, for $\theta\in {\cal B}_n$, if $u=\lambda_0(\theta)/\|\lambda_0(\theta)\|$, usual calculations (see \cite{Owen90} for instance) lead to 
\begeq{inegetape1}
 \|\lambda_0(\theta)\| \left( u^tB_n(\theta) u - M_n(\theta) \|G_n(\theta)\|\right) \leq \|G_n(\theta)\|
\fineq
for any $n$. Statements (\ref{Gntheta0TLC}), (\ref{controleMntheta0}), (\ref{cvgceunifBntheta}) and 
(\ref{controleMntheta}), (\ref{cvgceunifBntheta}), (\ref{controleGntheta}) thus respectively yield
\begeq{controlelambdatheta}
 \|\lambda_0(\theta_0)\| \; =\; O(k_n^{-1/2}) \makebox[2.cm][c]{a.s. \ \ and } \|\lambda_0(\theta)\| \; = \;  O(k_n^{-1/3}) \mbox{\ a.s. uniformly on ${\cal B}_n$}.
\fineq
Consequently, if we note $\gamma_{i,n}(\theta):= \ <\lambda_0(\theta),g_0(Z_{i,n},\theta)>$, we have 
\begeq{controlegammaintheta}
 \max_{i\leq n} |\gamma_{i,n}(\theta)| \leq \|\lambda_0(\theta)\|\, M_n(\theta) = o(1) \mbox{\ a.s. and uniformly on ${\cal B}_n$}
\fineq
and, using $(1+x)^{-1}=1-x+x^2(1+x)^{-1}$ and the identity $0=k_n^{-1}\sum_{i=1}^{k_n} (1+\gamma_{i,n}(\theta))^{-1} g_0(Z_{i,n},\theta)$,
we readily have 
\[
 \lambda_0(\theta) \; = \; (B_n(\theta))^{-1}G_n(\theta) + (B_n(\theta))^{-1}R_n(\theta)
\]
where  $R_n(\theta) = k_n^{-1} \sum_{i=1}^{k_n} (1+\gamma_{i,n}(\theta))^{-1}(\gamma_{i,n}(\theta))^2 g_0(Z_{i,n},\theta)$.
\ Since for $n$ sufficiently large we have $\|R_n(\theta)\| \leq 2k_n^{-1} \sum_{i=1}^{k_n} \|\lambda_0(\theta)\|^2\|g_0(Z_{i,n},\theta)\|^3
\leq 2\|\lambda_0(\theta)\|^2 M_n(\theta) \, tr(B_n(\theta))$, relations (\ref{controleMntheta0}), 
(\ref{controleMntheta}), (\ref{cvgceunifBntheta}) and (\ref{controlelambdatheta}) imply the following crucial relations
(the second one holding uniformly in $\theta\in{\cal B}_n$)
\begeq{devtlambdatheta}
 \lambda_0(\theta_0) = (B_n(\theta_0))^{-1} G_n(\theta_0) + o(k_n^{-1/2})
 \makebox[2.cm][c]{ a.s. \ and }
 \lambda_0(\theta) = (B_n(\theta))^{-1} G_n(\theta) + o(k_n^{-1/3}) \mbox{\ a.s.}.
\fineq
Using the Taylor expansion $\log(1+x)=x-\frac 1 2 x^2 + \frac 1 3 x^3(1+\xi)^{-3}$ (for some $\xi$ between $0$ and $x$) and statement (\ref{controlegammaintheta}),
we can proceed as above and obtain 
\[
 l_0(\theta) \; = \; 2\sum_{i=1}^{k_n} \log(1+\gamma_{i,n}(\theta)) \; = \; 
 2 \sum_{i=1}^{k_n} \gamma_{i,n}(\theta) - \sum_{i=1}^{k_n} (\gamma_{i,n}(\theta))^2 \, + \;  R'_n(\theta)
\]
where, for $n$ sufficiently large, $\|R'_n(\theta)\| \leq \frac {16} 3 \sum_{i=1}^{k_n} |\gamma_{i,n}(\theta)|^3 = 
o(1) \sum_{i=1}^{k_n} (\gamma_{i,n}(\theta))^2$. 
Using relation (\ref{devtlambdatheta}) as well as (\ref{cvgceunifBntheta}) and (\ref{controleGntheta}), 
we have for $\|\theta-\theta_0\|\leq k_n^{-1/3}$,
\begin{eqnarray*}
 \sum_{i=1}^{k_n} \gamma_{i,n}(\theta) 
    \; = \; k_n \lambda_0(\theta)^t G_n(\theta) & = & k_n(G_n(\theta)+o(k_n^{-1/3}))^t\, (B_n(\theta))^{-1}\, G_n(\theta) 
 \\ & = & k_n(G_n(\theta))^t\, B^{-1}\, G_n(\theta) + o(k_n^{1/3})
\end{eqnarray*}
and similarly, $\sum_{i=1}^{k_n} (\gamma_{i,n}(\theta))^2 \; = \; k_n(G_n(\theta))^t\, B^{-1}\, G_n(\theta) + o(k_n^{1/3})$. Therefore,
if $\theta=\theta_0+uk_n^{-1/3}$ with $\|u\|=1$, using (\ref{controleGntheta}) and the almost sure bound (\ref{Gntheta0vitesseps}) 
for $G_n(\theta_0)$, we obtain
\begin{eqnarray*}
 l_0(\theta) 
 & = & k_n(G_n(\theta))^t\, B^{-1}\, G_n(\theta) + o(k_n^{1/3}) \\
 & = & k_n^{1/3}[k_n^{1/3}G_n(\theta_0) + (A+o(1))u]^t\, B^{-1} [k_n^{1/3}G_n(\theta_0) + (A+o(1))u] + o(k_n^{1/3}) \\
 & = & k_n^{1/3} u^t A^tB^{-1}A u + o(k_n^{1/3}).
\end{eqnarray*}
Consequently, if $a>0$ denotes the smallest eigenvalue of $A^tB^{-1}A$ and $\epsilon\in]0,a[$, we have for $n$ sufficiently large  
\[
 l_0(\theta)\; \geq \; (a-\epsilon)k_n^{1/3} \mbox{\ almost surely and uniformly for $\theta\in cl({\cal B}_n)$.}
\]
On the other hand, we obtain in a similar manner
\[
 l_0(\theta_0) \; = \; (\sqrt{k_n}G_n(\theta_0))^t B^{-1} (\sqrt{k_n}G_n(\theta_0)) + o(1) \ \ \mbox{(a.s.)}
\]
which converges in distribution to $\chi^2(2)$ in view of (\ref{Gntheta0TLC}), and is also $o(\log\log k_n)$ almost surely, thanks to 
(\ref{Gntheta0vitesseps}). 
\zun\\
We have thus proved the following 

\begin{prop}
\label{proplocalminimizer}
Under the conditions of Theorem \ref{asymptprofile}, and as $n\tinf$ with probability one, the empirical log-likelihood ratio $l_0(\cdot)$ 
admits a local minimizer $\hat\theta$ in the interior of the ball ${\cal B}_n=\{ \, \theta\in\bR \, ; \, \|\theta-\theta_0\| \leq k_n^{-1/3} \, \}$.
This means that almost surely, for $n$ large, there exists a local minimizer $\hat\sigma_{\gamma_0}$ of the profile empirical log-likelihood 
$\sigma\mapsto l(\gamma_0,\sigma)$ such that $\hat\sigma_{\gamma_0}/\sigma_{0n}$ is close to $1$ with rate $k_n^{-1/3}$. 
\end{prop}

Now that we have identified some empirical likelihood estimator $\hat\theta$ and proved it consistently estimates $\theta_0$,
we want to identify its asymptotic distribution, which will enable us to obtain the convergence in distribution of $l_0(\hat\theta)$ towards 
$\chi^2(1)$. 
\zun\\
As it is done in \cite{QinLawless94}, we introduce the functions defined on $\bR\times\bR^2$  
\begin{eqnarray*}
 Q_{1,n}(\theta,\lambda) & = & \frac 1 {k_n} \sum_{i=1}^{k_n} (1+<\lambda,g_0(Z_{i,n},\theta)>)^{-1} g_0(Z_{i,n},\theta) \zun\\
 Q_{2,n}(\theta,\lambda) & = & \frac 1 {k_n} \sum_{i=1}^{k_n} (1+<\lambda,g_0(Z_{i,n},\theta)>)^{-1}  
                                  \left(\frac{\partial g_0}{\partial\theta}(Z_{i,n},\theta)\right)^t \lambda 
\end{eqnarray*}
and see that ($\forall\theta$) $Q_{1,n}(\theta,\lambda_0(\theta))=0$ (by definition of $\lambda_0(\theta)$), $Q_{1,n}(\theta,0)=G_n(\theta)$, and 
$Q_{2,n}(\theta,\lambda_0(\theta)) = (\partial l_0/\partial\theta)(\theta)$, which is null at $\theta=\hat\theta$.
\zun\\
A Taylor expansion of $Q_{1,n}$ and $Q_{2,n}$ between $(\theta_0,0)$ and $(\hat\theta,\lambda_0(\hat\theta))$ shows that
there exists some $(\theta^*_n,\lambda^*_n)$ satisfying $\|\theta^*_n-\theta_0\|\leq\|\hat\theta-\theta_0\|\leq k_n^{-1/3}$,
$\|\lambda^*_n\|\leq\|\lambda_0(\hat\theta)\| = O(k_n^{-1/3})$ (thanks to Proposition \ref{proplocalminimizer} and 
(\ref{controlelambdatheta})), and such that
\[
 \left( \begar{c} 0 \\ 0 \finar \right) 
 \; = \; 
 \left( \begar{c} Q_{1,n}(\hat\theta,\lambda_0(\hat\theta)) \\ Q_{2,n}(\hat\theta,\lambda_0(\hat\theta)) \finar \right) 
 \; = \; 
 \left( \begar{c} Q_{1,n}(\theta_0,0) \\ Q_{2,n}(\theta_0,0) \finar \right) \, - \, S_n(\theta^*_n,\lambda^*_n) 
 \left( \begar{c} \hat\theta-\theta_0 \\ \lambda_0(\hat\theta) \finar \right)
\]
where
\[
 S_n(\theta,\lambda) \; := \; \left. \left( \begar{cc} -\partial Q_{1,n}/\partial\theta & -\partial Q_{1,n}/\partial\lambda \\
  -\partial Q_{2,n}/\partial\theta & -\partial Q_{2,n}/\partial\lambda \finar \right) \right|_{(\theta,\lambda)}
\]
Differential calculus leads to 
\[ 
 S_n(\theta_0,0) \; = \; \left( \begar{cc} -\frac 1 {k_n} \sum_{i=1}^{k_n} \frac{\partial g}{\partial\theta}(Z_{i,n},\theta_0) 
 & \frac 1 {k_n} \sum_{i=1}^{k_n} g(Z_{i,n},\theta_0)(g(Z_{i,n},\theta_0))^t \zun\\
 0 & -\frac 1 {k_n} \sum_{i=1}^{k_n} \left(\frac{\partial g}{\partial\theta}(Z_{i,n},\theta_0)\right)^t \finar \right)
\] 
thus, defining $V:=(A^tB^{-1}A)^{-1}$, relations (\ref{cvgceunifBntheta}) and (\ref{cvgceunifAntheta}) imply that $S_n(\theta_0,0)$ 
converges to the matrix
\[
 S \; = \; \left( \begar{cc} -A & B \\ 0 & -A^t \finar \right) 
\]
which is invertible with 
\[ 
 S^{-1} \; = \; \left( \begar{cc} C & D \\ E & F \finar \right)
 \; := \; \left( \begar{cc} -VA^tB^{-1} & -V \\ B^{-1}(I-AVA^tB^{-1}) & -B^{-1}AV \finar \right)  
\]
After tedious calculations and use of many of the statements previously derived from Proposition \ref{propprelim},
it can be proved that $\|S_n(\theta^*_n,\lambda^*_n) - S_n(\theta_0,0)\| = o_{\bP}(1)$ as $n\tinf$. Consequently,  
we obtain, for $n$ sufficiently large
\begeq{presqueTLCthetachap}
 \left( \begar{c} \sqrt{k_n}(\hat\theta-\theta_0) \zun\\ \sqrt{k_n}\lambda_0(\hat\theta) \finar \right) 
 \; = \; 
 S^{-1} \left( \begar{c} \sqrt{k_n}G_n(\theta_0) + o_{\bP}(\sqrt{k_n}\delta_n) \\ o_{\bP}(\sqrt{k_n}\delta_n)  \finar \right)
 \makebox[1.6cm][c]{ where } \delta_n := \|\hat\theta-\theta_0\| + \|\lambda_0(\hat\theta)\|
\fineq
We already know that $\delta_n=O(k_n^{-1/3})$, but now (\ref{presqueTLCthetachap}) implies that $\delta_n=O(k_n^{-1/2})$ and therefore
we have proved that 
\begin{eqnarray}
 \sqrt{k_n}(\hat\theta-\theta_0) = C\sqrt{k_n}G_n(\theta_0) + o_{\bP}(1) & \stackrel{d}{\ra} & {\cal N}(0,CBC^t) = {\cal N}(0,(A^tB^{-1}A)^{-1})
 \label{TLCthetachap} \zun\\ 
 \sqrt{k_n}\lambda_0(\hat\theta) = E\sqrt{k_n}G_n(\theta_0) + o_{\bP}(1) & \stackrel{d}{\ra} & {\cal N}(0,EBE^t) = {\cal N}(0,E)
 \label{TLClambdachap}
\end{eqnarray}
where we have used the fact that the matrix $E$ is symmetric and such that $EBE=E$, because $I-AVA^tB^{-1}$ is idempotent (note that 
the rank of $E$ is $2$ minus that of $AVA^tB^{-1}$, {\it i.e.} $rank(E)=1$). 
\zun\\
Applying relation (\ref{controleGntheta}) to $\theta=\hat\theta$, relation (\ref{TLCthetachap}) yields
\begin{eqnarray*}
 \sqrt{k_n}G_n(\hat\theta) \; = \; \sqrt{k_n}G_n(\theta_0) + A\sqrt{k_n}(\hat\theta-\theta_0) + o_{\bP}(1) & = &
 (I+AC)\sqrt{k_n}G_n(\theta_0) + o_{\bP}(1) \\ 
 & = & BE(\sqrt{k_n}G_n(\theta_0) + o_{\bP}(1))
\end{eqnarray*}
and this leads to the following appropriate development for $l_0(\hat\theta)$, using (\ref{TLClambdachap}) and (\ref{cvgceunifBntheta}) : 
\begin{eqnarray*}
 l_0(\hat\theta) & = & 2k_n(\lambda_0(\hat\theta))^t G_n(\hat\theta) - k_n (\lambda_0(\hat\theta))^t B_n(\hat\theta)\lambda_0(\hat\theta) +   
                        R'_n(\hat\theta) \zun\\
                 & = & (\sqrt{k_n}G_n(\theta_0)+o_{\bP}(1))^t(EBE)(\sqrt{k_n}G_n(\theta_0)+o_{\bP}(1)) + o_{\bP}(1) + R'_n(\hat\theta) \zun\\
                 & = & (\sqrt{k_n}G_n(\theta_0))^tE(\sqrt{k_n}G_n(\theta_0)) + o_{\bP}(1) + R'_n(\hat\theta)                  
\end{eqnarray*}
where $|R'_n(\hat\theta)|\leq o_{\bP}(1) \sum_{i=1}^{k_n} (\gamma_{i,n}(\hat\theta))^2 =  
o_{\bP}(1) \left( (\sqrt{k_n}G_n(\theta_0))^tE(\sqrt{k_n}G_n(\theta_0)) + o_{\bP}(1) \right)$. 
\zun\\
According to proposition $(viii)$ p.524 of \cite{Rao84}, since $\sqrt{k_n}G_n(\theta_0)$ converges in distribution to ${\cal N}(0,B)$, 
and $EBE=E$ with $rank(EB)=1$, the quadratic form $(\sqrt{k_n}G_n(\theta_0))^tE(\sqrt{k_n}G_n(\theta_0))$ converges in distribution to 
$\chi^2(2-1)=\chi^2(1)$, and Theorem \ref{asymptprofile} is proved.
\zcinq
                        
\subsection{Proof of Proposition \ref{propprelim}}
\label{preuveprop}

Note that throughout the whole proof, we will write $\gamma$ instead of $\gamma_0$ for convenience.

\subsubsection{Proof of (\ref{Gntheta0vitesseps})  and (\ref{Gntheta0TLC})}

\noindent 
Let us define 
\[
Z_i=\frac{U_i^{-\gamma}-1}{\gamma} \mbox{ and }  \Delta_i(\theta)=g_0(Z_{i,n},\theta)-g_0(Z_i,\theta),
\]
so that $(Z_i)_{1\leq i\leq k_n}$ is an i.i.d. sequence with common distribution GPD$(\gamma, 1)$ and 
\begeq{Gntheta0}
 G_{n} (\theta_0)= \frac{1}{k_n} \sum_{i=1}^{k_n} g_0(Z_i,\theta_0) + \frac{1}{k_n} \sum_{i=1}^{k_n} \Delta_i(\theta_0).
\fineq
If $r < \frac{1}{2}$, $B$ is well defined as the covariance matrix of $g_0(Z_1,\theta_0)$ (a straightforward calculation gives the expression of $B$), 
and consequently the LIL and CLT  imply that 
\[
 \frac{1}{k_n}  \sum_{i=1}^{k_n} g_0(Z_i,\theta_0) =O\left( k_n^{-1/2}(\log\log k_n)^{1/2} \right)  \ a.s.  \ \ 
 \mbox{ and }  \ \ 
 \frac{1}{\sqrt{k_n}}  \sum_{i=1}^{k_n} g_0(Z_i,\theta_0) \stackrel{d}{\longrightarrow} \mathcal{N} \left(0, \, B \right). 
\]
Therefore, according to $(\ref{Gntheta0})$ and to the assumption $\sqrt{k_n} a_n \rightarrow 0$ (as $n \rightarrow +\infty$), 
in order to prove (\ref{Gntheta0vitesseps})  and (\ref{Gntheta0TLC}) it remains to establish that 
\begeq{controledeltamoyen}
 \frac{1}{k_n} \sum_{i=1}^{k_n} \Delta_i(\theta_0) =O(a_n) \ \ a.s.
\fineq
Since we can take $\sigma_{0n} := \sigma_0(u_n) = \frac{1}{\bar{F}(u_n)}  U^{\prime} \left( 1/ \bar{F}(u_n) \right)$, and 
recalling that we consider $Y_i=U(1/(U_i\bar{F}(u_n)))-U(1/\bar{F}(u_n))$, the application of the Potter-type bounds (\ref{bornesPotter}) 
to $t=1/\bar{F}(u_n)$ and $x=1/U_i$ yields, for all $\epsilon >0$ and $n$ sufficiently large, 
\begeq{Potter-type}
  (1-\epsilon) U_i^{\epsilon} K_{\gamma,\rho}(1/U_i) \leq \frac{Z_{i,n}-Z_i}{a_n}\leq (1+\epsilon) U_i^{-\epsilon} K_{\gamma,\rho}(1/U_i)  \ \ a.s.
\fineq
In the sequel, we will consider $a_n >0$  for large $n$ (the case $a_n <0$ being similar) and note $K_i=K_{\gamma,\rho}(1/U_i)$, as well
as $\Delta_i^1(\theta)$ and $\Delta_i^2(\theta)$ the two components of $\Delta_i(\theta)$. 
\zdeux

\begitem
\item[(i)] Control of $\Delta^1_i(\theta_0)$  
\[
  \Delta^1_i(\theta_0)  \; = \; \ln(1+\gamma Z_{i,n})-\ln(1+\gamma Z_i) \; = \; \ln \left( 1+\gamma U_i^{\gamma} \left(Z_{i,n} - Z_i  \right) \right).
\]
Use of $(\ref{Potter-type})$ leads to the following bounds (for all $\epsilon >0$ and $n$ sufficiently large),
\begeq{bornes1}
  \frac{1}{a_n} \ln \left(  1+a_n \gamma (1-\epsilon) U_i^{\gamma+\epsilon} K_i \right) \leq  \frac{\Delta^1_i (\theta_0)}{a_n} 
  \leq \frac{1}{a_n} \ln \left(  1+a_n \gamma (1+\epsilon) U_i^{\gamma-\epsilon} K_i \right) \ \  a.s. 
\fineq
If we set $B_i^+:=\gamma (1+\epsilon) U_i^{\gamma-\epsilon} K_i$ and $B_i^-:=\gamma (1-\epsilon) U_i^{\gamma+\epsilon} K_i$, 
Lemmas \ref{KiVi} and \ref{lemmecontrolemax} (stated and proved in the Appendix) imply that $B_i^+$ and $B_i^-$ are both square integrable and therefore 
$\max_{i\leq k_n} a_n B_i^+$ and $\max_{i\leq k_n} a_n B_i^-$ are both, almost surely, $o(\sqrt{k_n}a_n)$, which is $o(1)$ according 
to our assumption on $(k_n)$. 
\zun\\
Consequently, the inequality $\frac{2}{3} x \leq \ln(1+x) \leq x \ (\forall x\in [0,1/2])$ yields the following bounds, for all $\epsilon >0$ and 
$n$ sufficiently large, 
\begeq{bornes1bis}
\frac{2}{3} \frac{1}{k_n} \sum_{i=1}^{k_n} B_i^- \leq  \frac{1}{k_n} \sum_{i=1}^{k_n}  \frac{\Delta^1_i(\theta_0) }{a_n}  
 \leq \frac{1}{k_n} \sum_{i=1}^{k_n} B_i^+  \  \ a.s. 
\fineq
and therefore, for every $\epsilon>0$,
\begin{eqnarray*}
\frac{2}{3}  \gamma (1-\epsilon) \bE(U_1^{\gamma +\epsilon}K_1) \leq & \liminf \frac{1}{k_n} \sum_{i=1}^{k_n} \frac{\Delta^1_i(\theta_0) }{a_n} & 
\\  
\leq & \limsup \frac{1}{k_n} \sum_{i=1}^{k_n} \frac{\Delta^1_i(\theta_0) }{a_n} & \leq  \gamma (1+\epsilon) \bE(U_1^{\gamma -\epsilon}K_1).
\end{eqnarray*}
Letting $\epsilon$ go to $0$ gives $a_n^{-1}k_n^{-1}\sum_{i=1}^{k_n} \Delta^1_i(\theta_0)=O(1)$  \  a.s.
\zdeux

\item[(ii)] Control of $\Delta^2_i(\theta_0)$ 
\[
 \Delta^2_i(\theta_0) \; = \; (1+\gamma Z_{i,n})^{r/\gamma}-(1+\gamma Z_i)^{r/\gamma} 
 \; = \; U_i^{-r} \left( \left(1+\gamma U_i^{\gamma} \left(  Z_{i,n} - Z_i \right) \right)^{r/\gamma}-1 \right). 
\]
In the case $r<0$ (the case $r >0$ is similar), use of  $(\ref{Potter-type})$ yields the following bounds for all $\epsilon >0$ and $n$ large
\begin{eqnarray}
 & & \frac{U_i^{-r}}{a_n} \left( \left( 1+(1+\epsilon) \gamma a_n U_i^{\gamma-\epsilon} K_i \right)^{r/\gamma}-1 \right) \nonumber\\
 & & \leq \frac{\Delta^2_i(\theta_0) }{a_n}  
 \leq  \frac{U_i^{-r}}{a_n} \left( \left( 1+(1-\epsilon) \gamma a_n U_i^{\gamma+\epsilon} K_i \right)^{r/\gamma}-1 \right) 
 \label{bornes2} 
\end{eqnarray}
The inequality $\alpha x \leq (1+x)^{\alpha}-1 \leq \alpha  c x  \ (\forall  x\in [0,1/2]$, where $c= (\frac{3}{2})^{\alpha-1}>0$ and $\alpha=r/\gamma<0$)
yields, for all $\epsilon >0$ and $n$ sufficiently large, 
\begeq{bornes1bis}
 r(1-\epsilon) \frac{1}{k_n} \sum_{i=1}^{k_n} U_i^{\gamma-r+\epsilon} K_i  \leq  
 \frac{1}{k_n} \sum_{i=1}^{k_n}  \frac{\Delta^2_i (\theta_0)}{a_n} \leq 
 rc(1+\epsilon) \frac{1}{k_n} \sum_{i=1}^{k_n} U_i^{\gamma-r-\epsilon} K_i  \ \  a.s. 
\fineq
Once again, Lemma \ref{KiVi} ensures that $\bE(U_i^{\gamma -r \pm\epsilon} K_i)$ and $\bE((U_i^{\gamma -r  \pm\epsilon} K_i)^2)$ are finite (because $r< 1/2$), hence
for sufficiently small  $\epsilon >0$ 
\begin{eqnarray*}
 r (1-\epsilon) \bE(U_1^{\gamma -r+\epsilon}K_1) & \leq & \liminf \frac{1}{k_n} \sum_{i=1}^{k_n} \frac{\Delta^2_i }{a_n} \\
 & \leq & \limsup \frac{1}{k_n} \sum_{i=1}^{k_n} \frac{\Delta^2_i(\theta_0) }{a_n} \leq  r c (1+\epsilon) \bE(U_1^{\gamma-r -\epsilon}K_1) \ \ a.s .
\end{eqnarray*}
Letting $\epsilon$ go to $0$ yields $a_n^{-1}k_n^{-1}\sum_{i=1}^{k_n} \Delta^2_i(\theta_0)= O(1)$ a.s. and therefore 
(\ref{controledeltamoyen}) is proved.
\finit

\subsubsection{Proof of (\ref{controleMntheta0}) and (\ref{controleMntheta})}

With  $\Delta_i(\theta)$ and $Z_i$ being defined as previously, we have
\[  
  M_n (\theta) \; = \;  \max_{i \leq k_n} ||g_0(Z_{i,n},\theta)|| 
  \; \leq \; \max_{i \leq k_n} ||g_0(Z_i,\theta)|| +  \max_{i \leq k_n}  ||\Delta_i(\theta)||.
\]
Since the variables $g_0(Z_i,\theta_0)$ are square integrable, it comes (Lemma \ref{lemmecontrolemax})   
\[
 \max_{i \leq k_n} ||g_0(Z_i,\theta_0)|| = o (\sqrt{k_n}) \ \ a.s.
\]
On the other hand, part 1 of  Lemma \ref{majorantG1} implies that for $\theta$ in ${\cal{B}}_n$, 
$||g_0 (z,\theta)||^3 \leq G_1(z)$, for every $z\geq 0$ and $n$ sufficiently large. 
Since the variables $G_1(Z_i)$ are i.i.d. and integrable (part 4  of  Lemma \ref{majorantG1}), using Lemma \ref{lemmecontrolemax} we thus have 
\[
 \max_{i \leq k_n} ||g_0(Z_i,\theta)|| = o (k_n^{1/3}) \  \ a.s.
\]
We can now conclude the proof of (\ref{controleMntheta0}) and (\ref{controleMntheta}) by 
showing that, uniformly for $\theta$ in ${\cal{B}}_n$, $\max_{i \leq k_n}  |\Delta^j_i(\theta)|$ tends to $0$ almost surely 
for $j=1$ or $2$. Reminding that $\gamma Z_i=U_i^{-\gamma}-1$, we can show that
\[
 \Delta^1_i(\theta) \; = \;  \ln\left(1+\gamma \frac{Z_{i,n}}{\theta}\right)-\ln\left(1+\gamma \frac{Z_{i}}{\theta}\right) \; = \; 
 \ln \left\{ 1+\gamma U_i^{\gamma} \left( 1+(\theta-1)U_i^{\gamma} \right)^{-1} 
 \left(Z_{i,n} - Z_i \right) \right\}.
\]
Let $\delta>0$ and $\theta  \in ]1-\delta, 1+\delta[$. Using $(\ref{Potter-type})$, we have the following bounds (for all $\epsilon >0$ and $n$ sufficiently large),
\begeq{bornes1theta}
\frac{1}{a_n} \ln \left(  1+a_n \gamma \left(\frac{1-\epsilon}{1+\delta} \right) U_i^{\gamma+\epsilon} K_i \right) \leq  \frac{\Delta^1_i (\theta)}{a_n} \leq \frac{1}{a_n} \ln \left(  1+a_n \gamma \left(\frac{1+\epsilon}{1-\delta} \right) U_i^{\gamma-\epsilon} K_i \right) \ \  a.s. 
\fineq
where we supposed that $a_n>0$ (the other case is very similar). Proceeding as for the handling of $\Delta^1_i (\theta_0)$, and using $\frac 2 3 x \leq \ln(1+x) \leq x$ for $x\in [0,1/2]$, we obtain : for all $\delta>0$, $\theta  \in ]1-\delta, 1+\delta[$, $\epsilon >0$ and $n$ sufficiently large, 
\begeq{bornes1bistheta}
 \frac 2 3  \left(\frac{1-\epsilon}{1+\delta} \right) \gamma U_i^{\gamma+\epsilon} K_i \leq    
 \frac{\Delta^1_i(\theta) }{a_n} 
 \leq \left(\frac{1+\epsilon}{1-\delta} \right)  \gamma U_i^{\gamma-\epsilon} K_i  \ \  a.s. 
\fineq
Since $(\ref{bornes1bistheta})$ ensures that $\frac{\Delta^1_i(\theta) }{a_n} $ is of constant sign, for $n$ large enough we have 
\[
\sup_{\theta \in {\cal{B}}_n} \max_{i \leq k_n}  |\Delta^1_i(\theta)| \leq \sqrt{k_n} a_n \gamma \left(\frac{1+\epsilon}{1-\delta} \right) \frac{\max_{i \leq k_n} U_i^{\gamma-\epsilon} K_i}{\sqrt{k_n} }  \  \ a.s.  
\]
We conclude using Lemmas \ref{KiVi} and \ref{lemmecontrolemax} and assumption $\sqrt{k_n}a_n\ra 0$. The proof for $\Delta^2_i(\theta)$ is very similar. 

\subsubsection{Proof of (\ref{cvgceunifAntheta}) and (\ref{controleGntheta})}

Recall that $A_n(\theta)= \frac{1}{k_n} \sum_{i=1}^{k_n}  \frac{\partial g_0}{\partial\theta} (Z_{i,n},\theta)$ and let $A^*_n(\theta):= \frac{1}{k_n} \sum_{i=1}^{k_n}  \frac{\partial g_0}{\partial\theta} (Z_i,\theta)$, where the $Z_i$ were introduced previously.  We write 
\begeq{decompAnthetamoinsA}
A_n(\theta) - A= \left( A_n(\theta)- A^*_n(\theta) \right) \ +  \left( A^*_n(\theta)- A^*_n(\theta_0) \right) \ +   \left( A^*_n(\theta_0)- A \right) 
\fineq
and we will handle separately the three terms on the right hand side above. 
The third term goes to $0$ a.s. according to the strong law of large numbers (SLLN for short) and by definition of the $Z_i$ and $A$.  
The same is true (uniformly in $\theta$) for the second term, since part 3 of Lemma \ref{majorantG1} implies  
\[  
\sup_{\theta \in {\cal{B}}_n } \|A^*_n(\theta)- A^*_n(\theta_0)\|  \leq  
  \left(\sup_{\theta \in {\cal{B}}_n } || \theta-\theta_0 ||\right)  \frac{1}{k_n} \sum_{i=1}^{k_n}  G_3(Z_i).
\]
and the SLLN applies, thanks to part 4 of Lemma \ref{majorantG1}. 
It remains to study the first term of (\ref{decompAnthetamoinsA}) uniformly in $\theta$ in order to prove (\ref{cvgceunifAntheta}). We have
\[
 A_n(\theta)- A^*_n(\theta)  = \frac{1}{k_n} \sum_{i=1}^{k_n}  \tilde\Delta_i(\theta), 
\]
where the two components of $\tilde\Delta_i(\theta)$ are 
\begin{eqnarray*}
 \tilde\Delta^1_i(\theta) & = & -\theta^{-2}  \gamma Z_{i,n} (1+\gamma Z_{i,n} / \theta )^{-1}  +  
     \theta^{-2}  \gamma Z_i (1+\gamma Z_i / \theta)^{-1}  \zun\\ 
 \tilde\Delta^2_i(\theta) & = & -r \theta^{-2} Z_{i,n} (1+\gamma Z_{i,n}/\theta)^{r/\gamma-1}  +  
     r \theta^{-2} Z_i (1+\gamma  Z_i/\theta)^{r/\gamma-1} .  
\end{eqnarray*}
We shall give details for $\tilde\Delta^1_i(\theta)$ and the case $a_n >0$ (the case $a_n<0$ and the treatment  of 
$\tilde\Delta^2_i(\theta)$ can be handled very similarly). 
Let $\delta >0$, $\theta \in ]1-\delta,  1+\delta [$, and $V_i$ denote $(1+\gamma Z_i/\theta)^{-1}$. Use of the Potter-type bounds (\ref{Potter-type}) leads to the following bounds (for all $\epsilon >0$ and $n$ sufficiently large), 
\[
 V_i \left(  1+ a_n \theta^{-1}(1+\epsilon) \gamma U^{-\epsilon}_i  K_i V_i \right) ^{-1} 
 \leq \left(1+\gamma Z_{i,n} /\theta \right)^{-1} \leq 
 V_i \left(  1+ a_n \theta^{-1}(1-\epsilon) \gamma U^{\epsilon}_i  K_i V_i \right) ^{-1}  \ a.s. 
\]
After multiplication by $-\theta^{-2}\gamma Z_{i,n}$ and another use of (\ref{Potter-type}), we obtain
\begin{eqnarray*}
  & & -\theta^{-2}\gamma Z_i V_i \left\{ (1+ a_n B^-_i)^{-1} -1 \right\}
         -a_n \theta^{-2}(1+\epsilon) \gamma U^{-\epsilon}_i K_i V_i  (1+ a_n B^-_i)^{-1}    \zun \\
  & & \hspace*{0.3cm} \leq \ \tilde\Delta^1_i(\theta) \ \leq \
  -\theta^{-2}\gamma Z_i V_i \left\{ (1+ a_n B^+_i)^{-1} -1 \right\} 
         -a_n \theta^{-2}(1-\epsilon) \gamma U^{\epsilon}_i K_i V_i  (1+ a_n B^+_i)^{-1}    
\end{eqnarray*}
where $B^-_i= \theta^{-1}(1-\epsilon) \gamma U^{\epsilon}_i K_i V_i$ and $B^+_i=\theta^{-1}(1+\epsilon) \gamma U^{-\epsilon}_i K_i V_i$. 
Let us handle the upper bound first. 
We find easily that $(1-\delta)U_i^{\gamma}\leq V_i\leq (1+\delta)U_i^{\gamma}$, and therefore, by Lemmas 
\ref{KiVi} and \ref{lemmecontrolemax}, and assumption $\sqrt{k_n} a_n \rightarrow 0$,
\begeq{controleanBiplus}
 0 \; \leq \; \sup_{|\theta-1|\leq\delta} \max_{i\leq k_n} \left( a_n B^+_i \right)  
    \; \leq \; (1+\epsilon)(1+\delta)(1-\delta)^{-1} \gamma a_n \max_{i\leq k_n} U^{\gamma-\epsilon}_i K_i \; = \; o(1) \ \ a.s.
\fineq 
Consequently, using $(1+x)^{-1}-1=-x(1+x)^{-1}$, for $n$ sufficiently large and uniformly in $\theta \in ]1-\delta,1+\delta [$, we find (almost surely)
\begin{eqnarray*}
 & & \textstyle \frac 1 {k_n} \sum_{i=1}^{k_n} \frac{\tilde\Delta^1_i(\theta)}{a_n}  \\
 & & \textstyle \leq \ 
 \frac{(1+\epsilon)(1+\delta)^2}{(1-\delta)^3} \frac 1 {k_n} \sum_{i=1}^{k_n} \left((1-U_i^{\gamma})U_i^{\gamma-\epsilon}K_i \right) \ - \ 
  \frac{(1-\epsilon)(1-\delta)}{2(1+\delta)^2} \gamma \frac 1 {k_n} \sum_{i=1}^{k_n} \left(U_i^{\gamma+\epsilon}K_i \right) 
  \ = \ O(1) 
\end{eqnarray*}
The lower bound can be handled in the same way. Note that (\ref{controleGntheta}) is an immediate consequence of (\ref{cvgceunifAntheta}).

\subsubsection{Proof of (\ref{cvgceunifBntheta})}

Recall that $B_n(\theta)= k_n^{-1} \sum_{i=1}^{k_n}   g_0 (Z_{i,n},\theta) g_0 (Z_{i,n},\theta)^t $ 
and let 
\[
 B^*_n(\theta):= k_n^{-1} \sum_{i=1}^{k_n}   g_0 (Z_i,\theta) g_0 (Z_i,\theta)^t,
\] 
where the $Z_i$ were introduced previously.  We write 
\begeq{decompBnthetamoinsB}
B_n(\theta) - B= \left( B_n(\theta)- B^*_n(\theta) \right) \ +  \left( B^*_n(\theta)- B^*_n(\theta_0) \right) \ +   \left( B^*_n(\theta_0)- B \right) 
\fineq
The third term in the relation above goes to $0$ a.s. according to the SLLN and by definition of the $Z_i$ and $B$.  
Let us deal with the second term. For $\theta\in{\cal B}_n$, there exists some $\theta_n^*$ between $\theta_0$ and $\theta$ such that
(using parts 1 and 2 of Lemma \ref{majorantG1})
\begin{eqnarray*}
 \|B^*_n(\theta)- B^*_n(\theta_0)\| 
 & \leq & \|\theta-\theta_0\| . \frac{2}{k_n} \sum_{i=1}^{k_n} \left\| \frac{\partial g_0}{\partial \theta}(Z_i,\theta^*_n) \right\| \, \|g_0(Z_i,\theta^*_n)\| \zun\\
 & \leq & k_n^{-1/3} \max_{i\leq k_n} (G_1(Z_i))^{1/3} \frac{2}{k_n} \sum_{i=1}^{k_n}  G_2(Z_i).
\end{eqnarray*}
Therefore, combining part 4 of Lemma \ref{majorantG1}, Lemma \ref{lemmecontrolemax}, and the SLLN, we see that 
$\|B^*_n(\theta)- B^*_n(\theta_0)\|$ almost surely goes to $0$ as $n\tinf$, uniformly in $\theta\in{\cal B}_n$.  
\zdeux\\
It remains to study the first term of (\ref{decompBnthetamoinsB}) uniformly in $\theta$ in order to prove (\ref{cvgceunifBntheta}). We have
\begin{eqnarray*}
 B_n(\theta)- B^*_n(\theta)  
 & = & \frac{1}{k_n} \sum_{i=1}^{k_n} \Delta_i^t (\theta) g_0(Z_i,\theta) \; + \; \frac{1}{k_n} \sum_{i=1}^{k_n} \Delta_i (\theta)(g_0(Z_i,\theta))^t   
 \; + \; \frac{1}{k_n} \sum_{i=1}^{k_n}  \Delta_i (\theta) (\Delta_i (\theta) )^t  \\ 
 & = & \Gamma^t_{1,n}(\theta) + \Gamma_{1,n}(\theta) + \Gamma_{2,n}(\theta)
\end{eqnarray*}
with
\[
 \Delta_i (\theta)(g_0(Z_i,\theta))^t \; = \; \left( \begar{ll}   
 \Delta^1_i (\theta) (\ln(1+\gamma \frac{Z_i}{\theta})-\gamma) & \Delta^1_i (\theta) ((1+\gamma \frac{Z_i}{\theta})^{r/\gamma}- \frac{1}{1-r}) \zun \\
 \Delta^2_i (\theta) (\ln(1+\gamma \frac{Z_i}{\theta})-\gamma) & \Delta^2_i (\theta) ((1+\gamma \frac{Z_i}{\theta})^{r/\gamma}- \frac{1}{1-r})
 \finar \right)
\]
and 
\[
  \Delta_i (\theta) (\Delta_i (\theta) )^t \; = \; \left( \begar{cc}   
 (\Delta^1_i (\theta))^2  & \Delta^1_i (\theta) \ \Delta^2_i (\theta) \zun \\
 \Delta^1_i (\theta) \ \Delta^2_i (\theta) & (\Delta^2_i (\theta))^2
 \finar \right).
 \]
Considering the first element of the matrix $\Gamma_{1,n}(\theta)$, we have
 \[
 \textstyle
 a_n^{-1} \left| \frac{1}{k_n} \sum_{i=1}^{k_n} \Delta^1_i (\theta) \left(\ln\left(1+\gamma \frac{Z_i}{\theta}\right)- \gamma\right) \right| 
 \; \leq  \; \sqrt{\frac{1}{k_n} \sum_{i=1}^{k_n} \left( \frac{\Delta^1_i (\theta)}{a_n} \right)^2} \ 
               \sqrt{\frac{1}{k_n} \sum_{i=1}^{k_n} \left(\ln\left(1+\gamma \frac{Z_i}{\theta}\right)- \gamma\right)^2} 
 \]
and, applying the Cauchy-Schwarz inequality too for dealing with the other elements of $\Gamma_{1,n}(\theta)$ and $\Gamma_{2,n}(\theta)$, 
the convergence $B_n(\theta)- B^*_n(\theta)\ra 0$ (uniformly for $\theta\in{\cal B}_n$) will be proved as soon as we show that the means 
over $i=1$ to $k_n$ of each of the following quantities are almost surely bounded uniformly for $\theta\in{\cal B}_n$ :
\[
\textstyle
  \left(\ln\left(1+\gamma \frac{Z_i}{\theta}\right)- \gamma\right)^2, \ \left(\left(1+\gamma \frac{Z_i}{\theta}\right)^{r/\gamma}- \frac{1}{1-r}\right)^2, \
  \left( \frac{\Delta^1_i (\theta)}{a_n} \right)^2, \  \left( \frac{\Delta^2_i (\theta)}{a_n} \right)^2, \ 
  \frac{\Delta^1_i (\theta) \Delta^2_i (\theta)}{(a_n)^2}.
\]
Using Lemma \ref{majorantG1}, we see that the first two elements of this list are both uniformly bounded by 
$\frac{1}{k_n} \sum_{i=1}^{k_n} (G_1(Z_i))^{2/3}$ which converges almost surely. On the other hand, 
according to relation (\ref{bornes1bistheta}) and since $\bE(U_1^{2\gamma-2\epsilon} K_1^2)$ is finite (by Lemma \ref{KiVi}),  
$k_n^{-1}\sum_{i=1}^{k_n} (\Delta^1_i (\theta)/a_n)^2$ is uniformly almost surely bounded. Similarly the same is true for
$k_n^{-1}\sum_{i=1}^{k_n} (\Delta^2_i (\theta)/a_n)^2$, as well as for $k_n^{-1}\sum_{i=1}^{k_n} \Delta^1_i (\theta) \Delta^2_i (\theta)/a_n^2$, 
and the proof of (\ref{cvgceunifBntheta}) is over. 


\subsection{Proof of Theorem \ref{asymptCouple}}
\label{preuve1}

We proceed as in the start of the proof of Theorem \ref{asymptprofile}, and consider that the variables $Y_i$ are the variables $\bar{F}^{-1}_{u_n}(U_i)$ where $(U_i)_{i\geq 1}$ is an i.i.d. sequence of standard uniform variables. Recall that 
\[
 l(\gamma,\sigma) = 2\sum_{i=1}^{k_n} \log \left( 1+<\lambda(\gamma,\sigma),g(Y_i,\gamma,\sigma)> \right) .
\]
Defining  $Z_{i,n}=Y_i/\sigma_{0n}$, $\theta_0=(\gamma_0,1)$, $\theta=(\gamma,s)$, and
\[
 \tilde l(\theta) =  2\sum_{i=1}^{k_n} \log \left( 1+<\tilde\lambda(\theta),g(Z_{i,n},\theta)> \right)
\]
where $\tilde\lambda(\theta)$ is such that 
\[
  \frac{1}{k_n} \sum_{i=1}^{k_n} \left(  1+<\tilde\lambda(\theta),g(Z_{i,n},\theta)>\right)^{-1} g(Z_{i,n},\theta) = 0,
\] 
it comes that $l(\gamma_0,\sigma_{0n})=\tilde l(\theta_0)$ since $g(Z_{i,n},\gamma,s)=g(Y_i,\gamma,\sigma_{0n}s)$.
We thus need to prove that $\tilde l(\theta_0)$ converges in distribution to $\chi^2(2)$. Following a very classical outline in empirical likelihood theory, 
it is easy to prove that this convergence is guaranteed as soon as we have the following statements (as $n\tinf$)
\begin{eqnarray*}
 \frac 1 {k_n} \sum_{i=1}^{k_n} g(Z_{i,n},\theta_0) \; \stackrel{\bP}{\lra} \; 0, \hspace{0.3cm}
 \frac 1 {k_n} \sum_{i=1}^{k_n} g(Z_{i,n},\theta_0)(g(Z_{i,n},\theta_0))^t \; \stackrel{\bP}{\lra} \; B,   \\
 \frac 1 {\sqrt{k_n}} \sum_{i=1}^{k_n} g(Z_{i,n},\theta_0) \; \stackrel{d}{\lra} \; {\cal N}(0,B),  \hspace{0.3cm}
 \max_{i\leq k_n} \|g(Z_{i,n},\theta_0)\| = o_{\bP}(\sqrt{k_n})
\end{eqnarray*} 
However, these statements are included in Proposition \ref{propprelim} and therefore Theorem \ref{asymptCouple} is proved. 
\zun\\
{\it Note } : Proposition \ref{propprelim} was stated under the assumption that $r<1/3$, but in fact $r<1/2$ is sufficient in order to prove
all the results concerning $\theta_0$ only and not for $\theta$ in a neighborhood of it (and the covariance matrix $B$ is well defined and invertible as 
soon as $r<1/2$).
\zdeux

\section{Conclusion}

This work deals with the problem of finding confidence regions for the parameters of the approximating GPD distribution in the
classical POT framework, for general heavy tailed distributions ($\gamma>0$). It is shown that
the application of the empirical likelihood (EL) method to the estimating equations of \citet{Zhang07} yields confidence regions with 
improved coverage accuracy in comparison to the Wald-type confidence regions issued from the CLT for some  estimators
of the couple $(\gamma,\sigma)$ (including the maximum likelihood estimator). It is also observed that coverage accuracy is
not always as good as one would expect, which means that this subject (and the related one of EL calibration) would need to be
further investigated.
\zdeux

A profile EL based confidence interval for the tail index is also obtained, and its performance in terms of coverage probability has
been compared to that of the confidence interval (CI) described in \citet{LuPeng02} and \citet{PengQi06} (which is known to perform
better than the Wald-type CI based on the CLT for the Hill estimator). In some simulations, the interval of Lu, Peng and Qi shows 
better performance, but in others this performance is limited to a very short range of number $k_n$ of excesses : this instability with 
respect to $k_n$ is much less present for the CI based on Zhang's equations.
\zdeux

We shall finish this conclusion with two remarks. The first is that some of the methodology of the proof of the profile EL result
(inspired by \citet{QinLawless94}) could prove useful in other settings (Proposition \ref{propprelim} lists properties which yield 
convergence in distribution of empirical likelihood ratio when the observations form a triangular array). 
The second remark is that plug-in calibration could be an interesting subject to investigate for obtaining CI for $\gamma$, in particular
in order to shorten computation time. 
\zdeux

\section{Appendix}

\begin{lem}
\label{majorantG1}  
Let $\gamma>0$, $r<1/3$ and, for $\theta >0$, 
\[
 g_0(z,\theta) :=  \left( \begar{cc}   g^1_0(z,\theta)  \zun \\  g^2_0(z,\theta)  \finar \right) 
 \ = \  \left( \begar{cc}  \ln(1+\gamma z/\theta) - \gamma \zun \\  (1+\gamma z/\theta)^{r/\gamma} - (1-r)^{-1}   \finar \right).
\]
If we consider,  for some positive constants $c_1$, $c_1'$, $c_2$, $c_3$ depending on $r$ and $\gamma$,
\[ \begar{ccc}
 G_1(z) & = &  c_1 \left( \ln (1+ \gamma z) \, + \, (1+ \gamma z)^{r/\gamma} \, + \, c_1' \right)^3  \zun \\ 
  G_2(z) & = &  c_2 \left( z (1+\gamma z)^{-1} \, + \, z  (1+\gamma z)^{r/\gamma-1} \right) \zun \\ 
  G_3(z) & = &   c_3 \left( z (1+\gamma z)^{-1} \, + \, z  (1+\gamma z)^{r/\gamma-1} \, + \, z^2 (1+\gamma z)^{-2} \, + \, z^2  (1+\gamma z)^{r/\gamma-2} \right). 
\finar \]
then there exists a neighborhood  of $\theta_0=1$ such that for all $\theta$ in this neighborhood and $\forall z \geq 0$, 
\begenum
\item  $||g_0 (z,\theta)||^3 \leq G_1(z)$
\item $|| \frac{\partial g_0 }{\partial \theta}(z,\theta) || \leq G_2(z)$
\item $|| \frac{\partial^2 g_0}{\partial^2 \theta}(z,\theta) || \leq G_3(z)$
\item If $Z$ has distribution $GPD(\gamma,1)$, then $\bE(G_j(Z))$ is finite for each $j \in \{1,2,3\}$.
\finenum
\end{lem}

\noindent{\it Proof of Lemma \ref{majorantG1}} :  we shall first give details for part 1, since parts 2 and 3 can be treated similarly.  
\zun\\
We shall bound  from above $|g^1_0(z,\theta)|$ and $|g^2_0(z,\theta)|$ in the neighborhood  of $\theta_0=1$. For $\delta >0$ 
and $\theta \in [1-\delta, 1+\delta]$, we have, 
\[
 \ln \left(1+ \frac{\gamma z}{1+\delta}\right) -\gamma \leq g^1_0(z,\theta) \leq \ln \left(1+ \frac{\gamma z}{1-\delta}\right) -\gamma
\]
and if $r<0$ (the case $r>0$ is similar)
\[
\left(1+ \frac{\gamma z}{1-\delta}\right)^{r/\gamma} -\frac{1}{1-r} \leq g^2_0(z,\theta) \leq \left(1+ \frac{\gamma z}{1+\delta}\right)^{r/\gamma} -\frac{1}{1-r}.
\]
According to Lemma \ref{lemmajorations}, if we take $\delta< \frac{1}{3}$, we thus have (for some positive constant $c$)
\[
|g^1_0(z,\theta)| \leq \ln (1+ \gamma z) + \gamma + \ln(3/2), \mbox{ \ and \ } |g^2_0(z,\theta)| \leq c (1+ \gamma z)^{r/\gamma} + \frac{1}{1-r}.  
\] 
This concludes the proof of part 1. 
\zun\\
If $U$ is uniformly distributed on $[0,1]$, then it is easy to check that the expectations $\bE\,[\,(\ln(1+\gamma Z))^a(1+\gamma Z)^{rb/\gamma}\,]
= \gamma^a \bE\,[\,(-\ln U)^a U^{-rb}\,]$ are finite for every $a$ and $b$ in $\{0,1,2,3\}$ because we assumed that $r<1/3$. 
Therefore $\bE(G_1(Z))$ is finite. Similar simple integral calculus leads to the same conclusion for $\bE(G_2(Z))$ and $\bE(G_3(Z))$.
\zdeux


\begin{lem}
\label{KiVi} 
Let $\gamma >0$, $\alpha \in \bR$, $\beta \geq 0$ and $U$ a uniform $[0,1]$ random variable .  
\begitem 
\item[$(i)$] If  $1-\alpha-\gamma >0$,  then $\bE(U^{-\alpha} K_{\gamma,\rho}(\frac{1}{U})(-\ln U)^{\beta})$ is finite. 
\item[$(ii)$] If $1-\alpha-2\gamma>0$, then $\bE(U^{-\alpha} K^2_{\gamma,\rho}(\frac{1}{U})(-\ln U)^{\beta})$ is finite. 
\finit
\end{lem}

\noindent{\it Proof of Lemma $\ref{KiVi}$} : we have 
\[ \begar{llll}
K_{\gamma,\rho}(x) & = & \frac{1}{\rho} \left(  \frac{x^{\gamma+\rho}}{\gamma+\rho} - \frac{x^{\gamma}}{\gamma}\right) + \frac{1}{\gamma +\rho} & \mbox{ if } \gamma+\rho \neq 0                                                                        \mbox{ and } \rho \neq 0 \\
  & = & -\frac{1}{\gamma} \left( \ln x - \frac{x^{\gamma}-1}{\gamma} \right)   & \mbox{ if }  \gamma+\rho= 0 \mbox{ and } \rho \neq 0 \\
  & = & \frac{1}{\gamma^2} \left( x^{\gamma} (\gamma \ln x -1) +1 \right)   & \mbox{ if } \rho= 0.
\finar \]
We consider statement $(i)$ and provide details only for the case $\gamma+\rho \neq 0$ and $\rho \neq 0$ 
(all the other cases being handled the same way). 
A simple change in variables readily gives
\[ 
\begar{lll}
 \bE(U^{-\alpha} K_{\gamma,\rho}(\frac{1}{U})(-\ln U)^{\beta}) 
 & =  & \frac{1}{\rho(\gamma+\rho)}  \int_0^{+\infty} \exp \left((\alpha+\gamma+\rho-1)y \right) y^{\beta} \ dy \zun \\
  & & - \frac{1}{\rho \gamma}  \int_0^{+\infty} \exp \left((\alpha+\gamma-1)y \right) y^{\beta} \ dy \zun \\
  & & + \frac{1}{\gamma(\gamma+\rho)}  \int_0^{+\infty} \exp \left((\alpha-1)y \right) y^{\beta} \ dy. 
\finar \]
But $ \int_0^{+\infty} e^{uy} \ y^{\beta} \  dy $ being finite if and only if $u<0$, this concludes the proof of $(i)$ since $\gamma$ is positive and $\rho$ is negative, in this case. The proof of statement $(ii)$ involves the same arguments and is thus omitted.
\zdeux

\begin{lem}
\label{lemmajorations}
Let $\gamma >0$, $z>0$ and $\delta \in [0, 1/3]$. Then, 
\[  
 \frac{1}{2} <  \frac{ 1+ \frac{\gamma z}{1\pm\delta} }{ 1+\gamma z } <  \frac 3 2 
 \makebox[1.3cm][c]{and}  \left| \ln \left(1+ \frac{\gamma z}{1\pm\delta} \right) - \ln (1+\gamma z) \right| < \ln(3/2)  
\] 
\end{lem}

\noindent
{\it Proof of Lemma $\ref{lemmajorations}$} :  since 
\[Ê  
 \frac{1+ \frac{\gamma z}{1+\delta}}{1+\gamma z} = 1 - \frac{\delta }{1+\delta} \frac{\gamma z}{1+\gamma z} 
 \mbox{ \ and \ } 
 \frac{1+ \frac{\gamma z}{1-\delta}}{1+\gamma z} = 1 + \frac{\delta }{1-\delta} \frac{\gamma z}{1+\gamma z} 
\]
it is clear that the first ratio is between $\frac{1}{2}$ and $1$ and the second one between $1$ and  $\frac{3}{2}$. The second statement comes from
\[
 \ln \left(1+ \frac{\gamma z}{1+d} \right) - \ln (1+\gamma z) \ = \ \ln \left( \frac 1 {1+d} \left( 1 + \frac d { 1+\gamma z } \right) \right)
\]
which absolute value is bounded by $\ln(4/3)$ for $d=\delta$ and by $\ln(3/2)$ for $d=-\delta$, since $\delta$ is assumed to be in $[0,1/3]$. 
\zdeux

\begin{lem}
\label{lemmecontrolemax}
Let $(k_n)$ be an integer sequence such that $k_n\ra +\infty$. If $(Z_i)_{i\geq 1}$ is an i.i.d. sequence of non-negative random variables
such that $\bE(|Z_1|^p)$ for some $p>0$, then 
\[
\max_{i\leq k_n} |Z_i| = o(k_n^{1/p})  \mbox{\ almost surely, \ as $n\tinf$}. 
\]
\end{lem}

\noindent{\it Proof of Lemma \ref{lemmecontrolemax}} : mimicking the proof of Lemma 11.2 in \cite{OwenBook}, we find that
$\max_{i\leq n} |Z_i| = o(n^{1/p})$ almost surely and thus it is also true on the subsequence $(k_n)$, so the lemma is proved. 
\zdeux

\end{document}